\documentclass[12pt,reqno]{amsart}

\usepackage{amsmath,amsthm,amsfonts,amscd,amssymb,euscript}
\numberwithin{equation}{section} 
\newcommand{\dbar}{\ensuremath{\bar \partial}}

\newcommand{\ad}{\ensuremath{\bar \partial^{*}  }}
\newcommand{\C}{\ensuremath{{\mathbb C}}}
\newcommand{\R}{\ensuremath{{\mathbb R}}}
\newcommand{\N}{\ensuremath{{\mathbb N}}}

\newcommand{\smooth}{\ensuremath{C^{\infty}}}
\newcommand{\smoothc}{\ensuremath{C_0 ^{\infty}}}
\newcommand{\ra}{\ensuremath{C^{\omega}}}
\newcommand{\I}{\ensuremath{\mathcal I}}
\newcommand{\II}{\ensuremath{\mathfrak I}}

\newcommand{\JJ}{\ensuremath{\mathfrak J}}

\newcommand{\KK}{\ensuremath{\mathfrak K}}
\newcommand{\V}{\ensuremath{\mathcal V}}

\newcommand{\M}{\ensuremath{\mathcal M}}
\newcommand{\E}{\ensuremath{\mathcal E}}
\newcommand{\Hf}{\ensuremath{\mathcal H}}

\begin{document}
\title[Coherence and Other Properties of Sheaves in the Kohn Algorithm]{Coherence and Other Properties of Sheaves \\ in the Kohn Algorithm}
\author{Andreea C. Nicoara}

\address{Department of Mathematics, University of Pennsylvania, 209 South $33^{rd}$ St.,  Philadelphia, PA 19104}

\email{anicoara@math.upenn.edu}

\subjclass[2010]{Primary  	32W05; 35A27; Secondary 32C05.}


\keywords{Kohn algorithm, real-analytic functions, quasi-flasque sheaf, quasi-coherent sheaf}

\begin{abstract}
In the smooth case, we prove quasi-flasqueness for the sheaves of all subelliptic multipliers as well as at each of the steps of the Kohn algorithm on a pseudoconvex domain in $\C^n.$ We use techniques by Jean-Claude Tougeron to show that if the domain has a real-analytic defining function, the modified Kohn algorithm involving generating ideals and taking real radicals only in the ring of real-analytic germs yields quasi-coherent sheaves. This sharpens a result obtained by J. J. Kohn in 1979.
\end{abstract}

\maketitle

\tableofcontents

\section{Introduction}

\bigskip Looking at the subellipticity of the $\dbar$-Neumann problem on pseudoconvex domains in $\C^n$ led Joseph J. Kohn to define subelliptic multipliers in \cite{kohnacta} and prove these formed an ideal sheaf closed under a particular radical operation that he called the real radical. This same operation is known to real algebraic geometry community as the \L ojasiewicz radical. The sheaf defined by Kohn consists of germs of $\smooth$ functions subject to certain conditions, which arise from the analysis of the $\dbar$-Neumann problem. Kohn constructed a refinement of these conditions that yields an algorithmic procedure for verifying the subellipticity of the $\dbar$-Neumann problem on pseudoconvex domains. The procedure known nowadays as the Kohn algorithm produces an increasing chain of sheaves of ideals of multipliers, whose termination at the sheaf of all germs of $\smooth$ functions on a neighborhood of the boundary of the domain implies subellipticity of the $\dbar$-Neumann problem at every point in that neighborhood.

The present work seeks to elucidate the properties of the sheaf of all subelliptic multipliers as well as of the sheaves of multipliers that occur at each step of the Kohn algorithm. The investigation will be carried out at first in the most general case, namely when the boundary of the domain is smooth, and then stronger conditions will be proven in the more particular case when the boundary of the domain is real-analytic for a very natural modification of the Kohn algorithm.

Kohn's sheaf of ideals of subelliptic multipliers is defined in \cite{kohnacta} over the non-Noetherian ring $\smooth.$ Fortunately, it turns out that this sheaf is quasi-flasque, a notion introduced by Jean-Claude Tougeron in \cite{tougeronqf}, and thus possesses  the chief property of a coherent sheaf in the context where stalks may not be finitely generated, namely that the sections generating the ideal of multipliers for $(p,q)$ forms $I^q(x_0)$ at some $x_0 \in \overline\Omega$ also generate $I^q(x)$ for $x$ sufficiently close to $x_0.$ Herein lies the motivation for this work. Any proof of the Kohn Conjecture, i.e. the statement that the termination of the Kohn algorithm on $(p,q)$ forms is equivalent to finite order of contact of $q$-dimensional complex varieties with the boundary of a smooth pseudoconvex domain, must use this sheaf property in an essential fashion.

\medskip
\newtheorem{mainthm}{Theorem}[section]
\begin{mainthm}
\label{qfsheafthm}  Let $\Omega$ in $\C^n$ be a domain with $\smooth$ boundary $b \Omega.$ Let $\tilde U$ be any open subset of $b \Omega$ such that its closure $\overline{\tilde U}$ is compact as a subset of $\C^n.$ $b \Omega$ may be taken as $\tilde U$ if the domain is bounded. Let $\II^q$ be the sheaf of subelliptic multipliers for the $\dbar$-Neumann problem on $(p,q)$ forms defined on $\tilde U.$ The sheaf $\II^q$ satisfies the following two properties:
\begin{enumerate}
\item[(a)] $\II^q$ is quasi-flasque;
\item[(b)] Let $x_0$ be any point in $b\Omega.$ If sections $s_j \in \II^q(\tilde U)$ generate $\II^q(x_0)$ for $j \in J,$ $J$ an indexing set, then $s_j$ also generate $\II^q(x)$ for $x$ sufficiently close to $x_0.$ 
\end{enumerate}
\end{mainthm}

Compactness is needed here in order to pass from local to global subelliptic multipliers.

\medskip
\newtheorem{qfsteps}[mainthm]{Theorem}
\begin{qfsteps}
Let $\Omega$ in $\C^n$ be a pseudoconvex domain with $\smooth$ boundary $b \Omega.$ Let  $\II^q_k$ be the multiplier ideal sheaf given by the Kohn algorithm on $b \Omega$ at step $k,$ for each $k \geq 1.$ The sheaf $\II^q_k$ satisfies the following two properties:
\begin{enumerate}
\item[(a)] $\II^q_k$ is quasi-flasque;
\item[(b)] Let $x_0$ be any point in $b\Omega.$ If sections $s_j \in \II^q_k(b \Omega)$ generate $\II^q_k(x_0)$ for $j \in J,$ $J$ an indexing set, then $s_j$ also generate $\II^q_k(x)$ for $x$ sufficiently close to $x_0.$
\end{enumerate}
\label{quasiflasquesteps}
\end{qfsteps}

The first ideal in the Kohn algorithm on $(p,q)$ forms is the real radical of the ideal generated by the defining function of the domain $r$ and all Levi minors $\text{coeff}\{\partial r \wedge \dbar r \wedge (\partial \dbar r)^{n-q}\},$ and these Levi minors are subelliptic multipliers only if the domain is pseudoconvex. Therefore, the Kohn algorithm has no meaning if the domain fails to be pseudoconvex. On the other hand, unlike Theorem~\ref{qfsheafthm}, no compactness is needed in Theorem~\ref{quasiflasquesteps} to get the quasi-flasque property for the sheaves at all steps of the Kohn algorithm since apart from the real radical operation, the functions that appear in the Kohn algorithm are generated by global sections.

The Kohn algorithm involves differentiating, computing determinants, and taking real radicals starting from the defining function of the domain $r.$ When this defining function is real-analytic, differentiating and computing determinants yields other real-analytic functions, hence it is only natural to modify the Kohn algorithm so that instead of generating ideals and taking real radicals in the ring of smooth germs $\smooth,$ these operations take place in the Noetherian ring of real-analytic germs $\ra.$ We will denote by $\tilde \II^q$ the sheaf of real-analytic subelliptic multipliers for the $\dbar$-Neumann problem on $(p,q)$ forms and by $\tilde \II^q_k$ the sheaf of real-analytic subelliptic multipliers obtained at step $k$ of the modified Kohn algorithm. A result of Tougeron in  \cite{tougeronqf} allows us to show quasi-coherence for $\tilde \II^q$ and $\tilde \II^q_k$ via the arguments used to establish the quasi-flasque property of $\II^q$ and $\II^q_k:$

\medskip
\newtheorem{fintyp}[mainthm]{Theorem}
\begin{fintyp}
Let $\Omega$ in $\C^n$ be a domain with real-analytic boundary $b \Omega.$ Let $\tilde U$ be any open subset of $b \Omega$ such that $\tilde U$ is contained in a compact semianalytic subset $Y$ of $b \Omega.$ If $\Omega$ is bounded, $b \Omega$ itself may be taken as $\tilde U.$ The ideal sheaf $\tilde \II^q$ of real-analytic subelliptic multipliers for the $\dbar$-Neumann problem on $(p,q)$ forms defined on $\tilde U$ is coherent. Additionally, if $\Omega$ is pseudoconvex, the multiplier ideal sheaf $\tilde \II^q_k$ given by the modified Kohn algorithm on $\tilde U$ at step $k$ for each $k \geq 1$ is also coherent. In other words, $\tilde \II^q$ and $\tilde \II^q_k$ for all $k \geq 1$ are quasi-coherent sheaves.  \label{finitetype}
\end{fintyp}

Since Yum-Tong Siu's finiteness theorem from \cite{siunoet} is used in a fundamental way in the proof of Theorem~\ref{finitetype}, nothing is known about coherence for the sheaves $\tilde \II^q$ and $\tilde \II^q_k$ in the absence of some type of compactness.

In the real-analytic case, Kohn proved in \cite{kohnacta} that if an ideal of germs of real-analytic functions $I$ at some $0 \in \R^p$ is closed under the real radical operation, then there is a sequence of points $\{x_p\}_{p=1,2, \dots}$ converging to $0$ such that around each point $x_p,$ there exists a neighborhood $U_p$ satisfying that for all points $y$ in the intersection of $U_p$ with the variety corresponding to $I$ denoted  $\V(I),$ the ideal of all germs at $y$ of real-analytic functions vanishing on $\V(I)$ denoted $\I_y \V(I)$ is generated by elements of $I.$ Combining Theorem~\ref{finitetype} with the {\L}ojasiewicz Nullstellensatz for real-analytic germs, we will show Kohn's conclusion holds on an entire neighborhood around each point $x_0 \in b \Omega$ for $\tilde \II^q$ and $\tilde \II^q_k$ for all $k \geq 1.$

The paper is organized as follows: Section~\ref{subkohnalg} recalls the central definitions related to the subellipticity of the $\dbar$-Neumann problem and introduces the
Kohn algorithm. Section~\ref{qfsect} defines quasi-flasque sheaves and proves Theorems~\ref{qfsheafthm} and \ref{quasiflasquesteps}. The real-analytic case is treated in Section~\ref{rasect}, where Theorem~\ref{finitetype} is proven.

The author is grateful to Francesca Aquistapace, Vasile Brinzanescu, Fabrizio Broglia, and Stephen S. Shatz for a number of crucial conversations regarding sheaves.

\bigskip\bigskip

\section{The Kohn algorithm} \label{subkohnalg}

\bigskip

We shall very briefly cover here subelliptic multipliers for the $\dbar$-Neumann problem and the Kohn algorithm. See \cite{kohnacta} for full details or \cite{levidet} for a comprehensive outline. We start with Kohn's definition of a subelliptic multiplier:

\smallskip
\newtheorem{subellmult}{Definition}[section]
\begin{subellmult}
Let $\Omega$ be a domain in $\C^n$ and let $x_0 \in
\overline{\Omega}.$ A $\smooth$ function $f$ is called a
subelliptic multiplier at $x_0$ for the $\bar\partial$-Neumann
problem on $\Omega$ if there exist a neighborhood $U$ of $x_0$ and
constants $C, \epsilon > 0$ such that
\begin{equation}
||\, f \varphi\, ||_{\,U, \epsilon}^2 \leq C \, ( \, ||\,\dbar \,
\varphi\,||^2_{\, 0} + ||\, \ad  \varphi \,||^2_{\, 0} +
||\,\varphi \,||^2_{\, 0} \,) \label{subellest}
\end{equation}
for all $(p,q)$ forms $\varphi \in \smoothc (U \cap \overline{\Omega}) \cap Dom (\ad),$ where $||\, \cdot \, ||_{\, U,\epsilon}$ is the local Sobolev norm of
order $\epsilon$ on $U$ and $||\, \cdot \,||_{\, 0}$ is the $L^2$ norm. \label{subellmultdef}
\end{subellmult}

\smallskip\noindent Let $I^q (x_0)$ be the set of all subelliptic
multipliers at $x_0.$ The holomorphic part of the forms makes no difference when gauging subellipticity, so we have dropped $p$ from the notation.  We are only interested in looking at $b \Omega,$ the boundary of the domain $\Omega,$ as the $\dbar$-Neumann problem is automatically subelliptic inside the domain.

Definition~\ref{subellmultdef} yields a sheaf of germs of subelliptic multipliers that detects the subellipticity of the $\dbar$-Neumann problem since capturing a function $f$ inside $I^q(x_0)$ satisfying $f(x_0) \neq 0$ implies the existence of a subelliptic estimate at $x_0.$ The underlying ring is $\smooth,$ so at each point $I^q (x_0) \subset \smooth_{b \Omega} (x_0),$ the ring of germs of smooth functions on $b \Omega$ at $x_0 \in b \Omega.$ Let us denote by ${\mathfrak G}_{b \Omega}^{\smooth}$ the sheaf of germs of smooth functions defined on the boundary of the domain $b \Omega$ and by $\II^q$ the sheaf of subelliptic multipliers for the $\dbar$-Neumann problem on $(0,q)$ forms.

\medskip\noindent Theorem $1.21$ in \cite{kohnacta} sums up the properties of subelliptic multipliers Kohn proved:

\medskip
\newtheorem{subellcor}[subellmult]{Theorem}
\begin{subellcor}
If $\Omega$ is pseudoconvex with a $\smooth$ boundary and if $x_0
\in \overline{\Omega},$ then we have:
\begin{enumerate}
\item[(a)] $I^q (x_0)$ is an ideal.
\item[(b)] $I^q (x_0) = \sqrt[\R]{I^q (x_0)}.$
\item[(c)] If $r=0$ on $b \Omega,$ then $r \in I^q (x_0)$ and the
coefficients of $\partial r \wedge \dbar r \wedge (\partial \dbar
r)^{n-q}$ are in $I^q (x_0).$
\item[(d)] If $f_1, \dots, f_{n-q} \in I^q (x_0),$ then the
coefficients of $\partial f_1 \wedge \dots \wedge \partial f_j
\wedge \partial r \wedge \dbar r \wedge (\partial \dbar
r)^{n-q-j}$ are in $I^q (x_0),$ for $j \leq
n-q.$\label{subellpropcor}
\end{enumerate}
\end{subellcor}

\noindent Note that the real radical $\sqrt[\R]{I^q (x_0)}$ is computed in the ring of germs $\smooth_{b \Omega} (x_0),$ i.e. $\sqrt[\R]{I^q (x_0)}$ consists of all $f \in \smooth_{b \Omega} (x_0)$ such that there exist an open set $U_{x_0}$ containing $x_0,$ some $m \in \N^* ,$ and $g \in I^q(x_0)$ satisfying $|f(x)|^m \leq |g(x)|$ for all $x \in U_{x_0} \cap b \Omega.$

\medskip\noindent {\bf The Kohn Algorithm:}

\medskip\noindent {\bf Step 1}  $$I^q_1(x_0) = \sqrt[\R]{(\, r,\,
coeff\{\partial r \wedge \dbar r \wedge (\partial \dbar
r)^{n-q}\}\, )}$$  Here $( \, \cdot \, )$
denotes the ideal generated by the functions inside the
parentheses in the ring $\smooth_{b \Omega} (x_0),$ and $coeff\{\partial r \wedge \dbar r \wedge (\partial\dbar r)^{n-q}\}$ is the determinant of the Levi form for $q=1,$ while for $q>1$ it consists of minors of the Levi form.

\medskip\noindent {\bf Step (k+1)} $$I^q_{k+1} (x_0) = \sqrt[\R]{(\, I^q_k
(x_0),\, A^q_k (x_0)\, )},$$ where $$A^q_k (x_0)= coeff\{\partial
f_1 \wedge \dots \wedge \partial f_j \wedge \partial r \wedge
\dbar r \wedge (\partial \dbar r)^{n-q-j}\}$$ for $f_1, \dots, f_j
\in I^q_k (x_0)$ and $j \leq n-q.$

Theorem~\ref{subellpropcor} implies $I^q_k (x_0) \subset I^q (x_0)$ at each step $k.$ Furthermore, using $I^q_k
(x_0)$ to generate $I^q_{k+1} (x_0)$ ensures that at each point the algorithm produces an increasing chain of ideals
$$I^q_1 (x_0) \subset I^q_2 (x_0) \subset \cdots$$ in the ring of germs $\smooth_{b \Omega} (x_0).$ Let $\II^q_k$ be the subelliptic multipliers at step $k$ of the Kohn algorithm on $b \Omega,$ which is a sheaf of ideals that is a subsheaf of both $\II^q$ and of ${\mathfrak G}_{b \Omega}^{\smooth}.$ Note that the ideal sheaves $\II^q$ and $\II^q_k$ satisfy $\sqrt[\R]{\II^q}=\II^q$ and $\sqrt[\R]{\II^q_k}=\II^q_k$ for all $k \geq 1$ as the real radical is computed germwise. Additionally, the same property is true for the presheaves. In other words, for any open set $U \subset b\Omega,$ $\sqrt[\R]{I^q(U)}=I^q(U)$ and $\sqrt[\R]{I^q_k(U)}=I^q_k(U)$ for each $k \geq 1.$

We need some global notation for the varieties corresponding to
ideals of multipliers: $$\V_k^q = \V (\II^q_k),$$ where $$\V (\II^q_k)= supp \: \frac{ {\mathfrak G}_{b \Omega}^{\smooth}}{\II^q_k}.$$ In other words, $\V_k^q$ consists of all the points $x$ where the stalk of $\II^q_k$ at $x$ does not contain a unit.  It should be noted that for all $k \geq 1,$ $(\V_k^q, {\mathfrak G}_{b \Omega}^{\smooth}/ \II^q_k)$ are schemes, and if for some $k \geq 1$ and some open set $U \subset b \Omega,$ $\V_k^q \cap U = \emptyset,$ then the $\dbar$-Neumann problem is subelliptic for $(p,q)$ forms at every point of $U.$ $\V_k^q \cap U = \emptyset$ signifies that the Kohn algorithm is already finished by step $k$ on $U,$ a condition that is called Kohn finite ideal type for the neighborhood $U.$

\section{Quasi-flasque sheaves}
\label{qfsect}

Tougeron first defined quasi-flasque sheaves in \cite{tougeronqf}. That material is summarized in
section $V.6$ of \cite{tougeron}. \cite{hartshorne} is a good reference for elementary material on sheaves including the definition of a presheaf. We will use Tougeron's machinery to prove Theorems~\ref{qfsheafthm} and \ref{quasiflasquesteps} in this section.

\medskip
\newtheorem{qfdef}{Definition}[section]
\begin{qfdef}
Let $\tilde U$ be a given non-empty open set in $\R^m,$ and let
$\E$ be the sheaf of $\smooth$ germs on $\tilde U.$ A sheaf $\M$
of $\E$-modules is called quasi-flasque if for every open set $U
\subset \tilde U$ the canonical homomorphism
$$\M (\tilde U) \otimes_{\E(\tilde U)} \E (U) \longrightarrow \M
(U)$$ is an isomorphism.
\end{qfdef}

\noindent {\bf Remarks:} 

\smallskip \noindent (1) $\M (\tilde U) = \Gamma (\tilde U, \M),$ the sections of $\M$ on $\tilde U,$ and similarly $\M (U) = \Gamma (U, \M).$ The same holds for $\E,$ i.e. $\E(\tilde U)=\Gamma (\tilde U, \E)=\smooth (\tilde U)$ and $\E(U)=\Gamma (U, \E)=\smooth (U).$

\smallskip \noindent (2) Essentially, the quasi-flasque property has to be proven on the presheaf.

\smallskip \noindent (3) The tensor product $\M (\tilde U) \otimes_{\E(\tilde U)} \E (U)$ can be more easily characterized in this context because the ring $\E (U)$ is a flat module over $\E(\tilde U).$

\medskip Before we can examine the Kohn multiplier ideal sheaves and prove they are quasi-flasque, we have to recall some algebraic results on the ring $\smooth$ concerning flatness. If $U$ and $\tilde U$ are two open sets such that $U \subset \tilde U,$ then $\smooth(U)$ is a flat module over $\smooth(\tilde U).$ Here flatness has its
algebraic meaning that Serre introduced in \cite{serre} and has nothing to do with the concept of flatness
defined for functions in analysis. In \cite{tougeron}
Tougeron obtains this flatness property of the ring of
smooth functions over the smaller neighborhood as a module over
the ring of smooth functions over the large neighborhood via the following technical lemma, which may be found in subsection $V.6$ on pages 113-4:

\medskip
\newtheorem{flatnesslemma}[qfdef]{Lemma}
\begin{flatnesslemma}
Let $U \subset \tilde U$ be open, and let $\{f_i \}_{i \in \N}$ be a
countable family of functions in $\smooth (U).$ There exists a
function $\alpha \in \smooth (\tilde U)$ such that \label{flatness}
\begin{enumerate}
\item $\alpha \equiv 0$ on $\tilde U - U$ and $\alpha (x) \neq 0$ for all
$x \in U;$
\item The functions $\alpha \cdot f_i$ can be extended as smooth
functions $f_i'$ in $\smooth (\tilde U)$ such that $$f_i' = \alpha \cdot
f_i \equiv 0$$ on $\tilde U-U.$
\end{enumerate}
\end{flatnesslemma}

\smallskip\noindent {\bf Remarks:} 

\smallskip\noindent (1) When we work with open sets $\tilde U$ such that $\overline{\tilde U}$ is compact in $\C^n,$ we can assume without loss of generality that $\alpha$ extends smoothly up to the boundary of $\tilde U.$ The reason is that $\overline{\tilde U}$ being compact in $\C^n$ implies there exists an open set $U'$ such that $\overline{\tilde U} \subset U'.$ By applying the lemma twice, first to the pair of sets $U$ and $\tilde U$ and a second time to the pair $\tilde U$ and $U'$ with the function $\alpha$ from the first application replacing the family $\{f_i \}_{i \in \N},$ it follows that $\alpha \in \smooth (\overline{\tilde U}).$

\smallskip\noindent (2) The lemma does not specify the sign of $\alpha,$ but it is clear we can take $\alpha (x) \geq 0$ for all $x \in \tilde U.$

\bigskip For the reader's convenience, we recall here one of the standard ways of checking flatness. Consider a ring $A$ and a module $E$ over $A.$ Let $f=(f_1, \dots, f_n) \in A^n.$ Following Malgrange in \cite{malgrange}, we denote by $R(f,E)$ the relations of $f$ in $E,$ i.e. the submodule of $E^n$ consisting of the $n$-tuples $(e_1, \dots, e_n)$ verifying $\sum_{i=1}^n f_i \, e_i =0.$ Similarly, $R(f,A)$ are all the relations of $f$ in $A.$ The following result is Proposition 4.2 on p.42 of \cite{malgrange}:

\medskip
\newtheorem{flatnessprop}[qfdef]{Proposition}
\begin{flatnessprop}
A module $E$ is flat over a ring $A$ iff for every $n$ and every $f \in A^n,$ we have $R(f,E)=R(f,A) \, E.$
\label{flatnessproposition}
\end{flatnessprop}

\medskip
\newtheorem{flatnesscor}[qfdef]{Corollary (Flatness Property)}
\begin{flatnesscor}
\label{flatnesscorollary} Let $U \subset \tilde U$ be open, then
$\smooth(U)$ is a flat module over $\smooth(\tilde U).$
\end{flatnesscor}

\smallskip\noindent {\bf Proof:} To show $\smooth(U)$ is a flat module
over $\smooth(\tilde U),$ we need to check the criterion from Proposition~\ref{flatnessproposition}. Let $\phi_1, \dots, \phi_t \in \smooth(\tilde U)$. We
consider any relation among the restrictions of these functions to
$\smooth(U),$ namely let $f_1, \dots, f_t \in \smooth (U)$ be such
that $$f_1 \, \phi_1 + \cdots + f_t \, \phi_t =0.$$ Apply Lemma~\ref{flatness} with $f_1, \dots, f_t$
replacing the countable family $\{f_i \}_{i \in \N},$ and consider
$f_i' \in \smooth (\tilde U)$ defined by $f_i' = \alpha \cdot f_i.$
$$f_1' \, \phi_1 + \cdots + f_t' \, \phi_t =0$$ is obviously a
relation on $\smooth (\tilde U).$ Since $\alpha$ is a unit on
$\smooth (U)$ by Lemma~\ref{flatness}, restrict the relation $f_1' \, \phi_1 + \cdots + f_t' \, \phi_t =0$ to $U$ and multiply by $\frac{1}{\alpha}$ to obtain the relation with which we started. Therefore, all relations on $\smooth(U)$ are indeed generated by relations on $\smooth(\tilde U)$ as needed. \qed

\medskip By Corollary~\ref{flatnesscorollary} for any $U \subset \tilde U$ open, $\E (U)$ is a flat $\E (\tilde
U)$-module, which implies $$\II (U) \otimes_{\E(\tilde U)} \E (U) \simeq \Gamma (\tilde U, \II) \, \E (U),$$ i.e. the ideal generated in the ring $\E(U)$ by all the sections in $\Gamma (\tilde U, \II).$ At times, it is more convenient to use the notation $\II_G (U) = \II(U) \otimes_{\E(\tilde U)} \E (U)$ for this object, where the subscript indicates we are dealing with an object generated by global sections. The collection $\II (U)=\Gamma(U, \II)$ for all open subsets $U$ of $\tilde U$ with the obvious maps represents the presheaf of $\II.$

Let us now examine estimate \eqref{subellest} from the definition of a subelliptic multiplier. We immediately see from this flatness construction that the sheaf $\II^q$ of all subelliptic multipliers has to be quasi-flasque:

\medskip
\newtheorem{bringhome}[qfdef]{Proposition}
\begin{bringhome}
\label{bringhomeall}
Let $\Omega$ in $\C^n$ be a domain with $\smooth$ boundary $b \Omega.$ Let $\tilde U$ be any open subset of $b \Omega$ such that its closure $\overline{\tilde U}$ is compact as a subset of $\C^n.$ $\tilde U$ may be taken to equal $b \Omega$ if the domain is bounded. Let $\II^q$ be the sheaf of subelliptic multipliers for the $\dbar$-Neumann problem on $(p,q)$ forms defined on $\tilde U.$ The sheaf $\II^q$ is quasi-flasque.
\end{bringhome}

\smallskip\noindent {\bf Proof:} Trivially, $\II^q_G (U) \subset \Gamma (U, \II^q),$ so quasi-flasqueness amounts to the opposite inclusion. Let $f \in \Gamma (U, \II^q).$ This means $f$ is a subelliptic multiplier on a neighborhood $U \ni x_0$ satisfying estimate \eqref{subellest}. As explained on p.93 of \cite{kohnacta}, without loss of generality we can employ the local tangential Sobolev norm on open subsets of the boundary of the domain $b \Omega$ instead of the local Sobolev norm on open subsets of $\C^n$ on the left-hand side of estimate \eqref{subellest}. Apply Lemma~\ref{flatness}, and consider
$\tilde f \in \smooth (\tilde U)$ defined by $\tilde f = \alpha \cdot f.$ Clearly,
\begin{equation*}
\frac{1}{M^2} \, ||\, \tilde f \varphi\, ||_{\,\tilde U, \epsilon}^2 \leq||\, f \varphi\, ||_{\,U, \epsilon}^2 \leq C \, ( \, ||\,\dbar \,
\varphi\,||^2_{\, 0} + ||\, \ad  \varphi \,||^2_{\, 0} +
||\,\varphi \,||^2_{\, 0} \,), 
\end{equation*}
where $M = \sup_{y \in \overline{\tilde U}} \, |\alpha (y)|.$ Since $\overline{\tilde U}$ is compact, $M < \infty.$ Therefore, $f$ is coming from a global section, i.e. $f \in \II^q_G (U)$ as needed.  \qed

\medskip\noindent In order to finish the proof of Theorem~\ref{qfsheafthm}, we recall  Tougeron's
Proposition 6.4 from section V.6 of his book \cite{tougeron} that if the stalk of a quasi-flasque sheaf is defined by global sections at a point $x_0,$ then those same sections generate the stalks of the sheaf at all neighboring points close enough to $x_0:$

\medskip
\newtheorem{jollygroup}[qfdef]{Proposition}
\begin{jollygroup}
\label{jollygroupprop}
Let $\M$ be a quasi-flasque sheaf on $\tilde U$ and let $a$ be a
point in $\tilde U.$ If sections $s_j \in \M (\tilde U)$ for $j
\in J$ generate $\M_a$, then these also generate $\M_x$ for $x$
sufficiently close to $a.$
\end{jollygroup}

\bigskip\noindent {\bf Proof of Theorem~\ref{qfsheafthm}} 

Part (a) follows from Proposition~\ref{bringhomeall}. 

Part (b) follows from Proposition~\ref{jollygroupprop}. \qed

\medskip A quick glance at the Kohn algorithm shows that except for the real radical operation it can be restated in terms of globally defined objects. We will then prove a more general lemma that starting with the presheaf of a quasi-flasque ideal sheaf $\II$ and closing its presheaf $\{ \II(U)\}_{U \subset \tilde U}$ via the real radical yields another quasi-flasque sheaf $\KK:$

\medskip
\newtheorem{qflemma}[qfdef]{Lemma}
\begin{qflemma}
Let $\tilde U$ be such that $\overline{\tilde U}$ is compact in $\R^m,$ and let $\II$ be an ideal sheaf on $\tilde U$ that is quasi-flasque. Consider the presheaf $\{ \II(U)\}_{U \subset \tilde U}$ of $\II$ and define the presheaf $\{ \KK(U)\}_{U \subset \tilde U},$ where $\KK (U)  =\sqrt[\R]{\II(U)}$ on any open $U \subset \tilde U.$ The sheaf $\KK$ corresponding to the presheaf $\{ \KK(U)\}_{U \subset \tilde U}$ is also quasi-flasque. \label{realradok}
\end{qflemma}

\smallskip\noindent {\bf Proof:} Since $\II$ and $\KK$ are subsheaves of $\E,$ our previous discussion on flatness means $\KK_G (U) \subset \Gamma (U, \KK)$ is trivially true and only the opposite inclusion needs to be established.  Let us then start with any $f \in \Gamma (U, \KK).$ We know $f \in \sqrt[\R]{\II(U)}$ by definition, so there exist some $m \in \N^* $ and $g \in \II(U)$ satisfying $|f(x)|^m \leq |g(x)|$ for all $x \in U.$ Since $\II$ is quasi-flasque, $g \in \II_G (U),$ i.e. there exist $g_1, \dots, g_t \in \II(\tilde U)$ and $\phi_1, \dots, \phi_t \in \smooth(U)$ such that $$g(x) = \phi_1(x) \, g_1(x) + \dots + \phi_t (x) \, g_t (x).$$ We now apply Lemma~\ref{flatness} with $f, \phi_1, \dots, \phi_t$
replacing the countable family $\{f_i \}_{i \in \N},$ and consider
$f', \phi'_1, \dots, \phi'_t \in \smooth (\tilde U)$ defined by $f' = \alpha \cdot f$ and $\phi_i' = \alpha \cdot \phi_i$ for $1 \leq i \leq t.$ Multiplying $|f(x)|^m \leq |g(x)|$ by $|\alpha(x)|^m$ and rearranging terms gives us $$|\alpha(x) \cdot f(x)|^m \leq |\alpha^{m-1}(x)\phi'_1(x) \, g_1 (x) + \dots + \alpha^{m-1}(x) \phi'_t(x) \, g_t (x)|$$ for all $x \in \tilde U$ as $\alpha(x) \geq 0$ everywhere on $\tilde U$ by the second remark following Lemma~\ref{flatness}. In other words, $\alpha \cdot f \in \sqrt[\R]{\II(\tilde U)},$ but $\sqrt[\R]{\II(\tilde U)}=\KK(\tilde U),$ so $\alpha \cdot f \in \KK(\tilde U).$ As $\alpha(x) \neq 0$ for all $x \in U,$ it follows that $f \in \KK_G(U)$ as needed. \qed

\medskip\noindent Let us now prove quasi-flasqueness for the sheaves at all steps of the Kohn algorithm:

\medskip
\newtheorem{bringhomek}[qfdef]{Proposition}
\begin{bringhomek}
\label{bringhomestepk}
Let $\Omega$ in $\C^n$ be a pseudoconvex domain with $\smooth$ boundary $b \Omega.$ Let  $\II^q_k$ be the multiplier ideal sheaf given by the Kohn algorithm on $b \Omega$ at step $k,$ for each $k \geq 1.$ The sheaf $\II^q_k$ is quasi-flasque for all $k \geq 1.$
\end{bringhomek}

\smallskip\noindent {\bf Proof:} We prove this statement by induction on the step $k.$

\smallskip\noindent {\bf Base case:} Germwise, the sheaf $\II^q_1$ is generated by $$I^q_1(x) = \sqrt[\R]{(\, r,\,
coeff\{\partial r \wedge \dbar r \wedge (\partial \dbar
r)^{n-q}\}\, )}.$$ Consider the sheaf $\JJ^q_1$ generated germwise by $$J^q_1(x)=(\, r,\,
coeff\{\partial r \wedge \dbar r \wedge (\partial \dbar
r)^{n-q}\}\, ).$$ The defining function of the domain $r$ and all the Levi minors $coeff\{\partial r \wedge \dbar r \wedge (\partial \dbar r)^{n-q}\}$ are global functions, so the sheaf $\JJ^q_1$ is trivially quasi-flasque. Apply Lemma~\ref{realradok} to conclude $\II^q_1$ is quasi-flasque as well.

\smallskip\noindent {\bf Inductive step:} Assume $\II^q_k$ is quasi-flasque. Germwise, the sheaf $\II^q_{k+1}$ is generated by $$I^q_{k+1} (x) = \sqrt[\R]{(\, I^q_k (x),\, A^q_k (x)\, )}.$$ Consider the sheaf $\JJ^q_{k+1}$ generated germwise by $$J^q_1(x)=(\, I^q_k (x),\, A^q_k (x)\, ).$$ By the inductive hypothesis, any element of $I^q_k(x)$ is generated by global sections. Now, $$A^q_k (x)= coeff\{\partial f_1 \wedge \dots \wedge \partial f_j \wedge \partial r \wedge \dbar r \wedge (\partial \dbar r)^{n-q-j}\}$$ for $f_1, \dots, f_j \in I^q_k (x)$ and $j \leq n-q.$ Once again, by the inductive hypothesis, each $f_l$ among $f_1, \dots, f_j$ is defined by global sections, i.e. there exist global sections $g_1^{(l)}, \dots, g_{i_l}^{(l)} \in \II^q_k(b \Omega)$ and smooth functions $a_1^{(l)}, \dots, a_{i_l}^{(l)} \in \smooth(U)$ such that $$f_l = a_1^{(l)} \, g_1^{(l)} + \dots + a_{i_l}^{(l)} \, g_{i_l}^{(l)}.$$ When computing the complex gradient $\partial f_l,$ the differentiation may fall on $a_1^{(l)}, \dots, a_{i_l}^{(l)},$ which will yield other functions still in $\smooth (U),$ or it may fall on the global sections $g_1^{(l)}, \dots, g_{i_l}^{(l)},$ which will yield other global sections in $ \smooth(b \Omega).$ Altogether, the elements of $A^q_k (x)$ are also coming from global sections, which makes the sheaf $\JJ^q_{k+1}$ quasi-flasque. Apply Lemma~\ref{realradok} to conclude $\II^q_{k+1}$ is also quasi-flasque. \qed

\bigskip\noindent {\bf Proof of Theorem~\ref{quasiflasquesteps}} 

Part (a) follows from Proposition~\ref{bringhomestepk}.

Part (b) follows from Proposition~\ref{jollygroupprop}. \qed

\bigskip\bigskip

\section{The real-analytic setting}
\label{rasect}

\bigskip

We now turn our attention to the case where the domain is defined by a real-analytic function $r.$ Following Kohn in section 6 of \cite{kohnacta}, we modify the algorithm so that ideals are generated and the real radical takes place only in the ring of real-analytic germs $\ra(x_0):$  $$I^q_1(x_0) = \sqrt[\R]{(\, r,\,coeff\{\partial r \wedge \dbar r \wedge (\partial \dbar r)^{n-q}\}\, )_{\ra(x_0)}}$$ and $$I^q_{k+1} (x_0) = \sqrt[\R]{(\, I^q_k (x_0),\, A^q_k (x_0)\, )_{\ra(x_0)}},$$ where $$A^q_k (x_0)= coeff\{\partial
f_1 \wedge \dots \wedge \partial f_j \wedge \partial r \wedge
\dbar r \wedge (\partial \dbar r)^{n-q-j}\}$$ for $f_1, \dots, f_j
\in I^q_k (x_0) \cap \ra(x_0)$ and $j \leq n-q.$ We then denote by $\tilde \II^q_k$ the sheaf of real-analytic subelliptic multipliers obtained at step $k$ of the modified Kohn algorithm for all $k \geq 1.$ Likewise, we denote by $\tilde \II^q$ the sheaf of real-analytic subelliptic multipliers for the $\dbar$-Neumann problem on $(p,q)$ forms, where this sheaf is generated by allowing only real-analytic multipliers in estimate \eqref{subellest}. Due to the fact mentioned in the proof of Proposition~\ref{bringhomeall} that the local tangential Sobolev norm can be substituted for the local Sobolev norm on the left-hand side of estimate \eqref{subellest}, we may view $\tilde \II^q$ and $\tilde \II^q_k$ for all $k \geq 1$ as real-analytic ideal sheaves on $b \Omega,$ a real-analytic manifold given by $\{ r = 0\},$ which is thus countable at infinity.

Tougeron explored in \cite{tougeronqf} when such real-analytic ideal sheaves defined on a real-analytic manifold countable at infinity are coherent. The canonical injection $\ra \rightarrow \E$ makes $\E$ a faithfully flat $\ra$-module. The reader may consult Malgrange's book \cite{malgrange} for additional information on faithful flatness as opposed to flatness, which we covered in Section~\ref{qfsect}. Furthermore, the functor $\Hf \leadsto \tilde \Hf = \Hf \otimes_{\ra} \E$ from the category of real-analytic sheaves on our real-analytic manifold countable at infinity to the category of sheaves of smooth germs on the same manifold is faithful and exact. Proposition VI on p.2973 of \cite{tougeronqf} gives the following criterion for coherence of $\Hf:$

\smallskip
\newtheorem{qfcoh}{Proposition}[section]
\begin{qfcoh}
A real-analytic sheaf $\Hf$ is coherent iff $\tilde \Hf$ is quasi-flasque and of finite type.
\label{qfcoherent}
\end{qfcoh}

\smallskip While the ring $\ra(x_0)$ of germs of real-analytic functions at a point $x_0 \in b \Omega$ is Noetherian, obviously the ring $\ra(U)$ of real-analytic functions on an open set $U \subset b \Omega$ is in general not Noetherian. To obtain coherence for the sheaves $\tilde \II^q$ and $\tilde \II^q_k$ for all $k \geq 1,$ some finiteness theorems are necessary. The following theorem was obtained by Siu in \cite{siunoet}. The reader is also directed to \cite{zame} for a different proof and to \cite{adkinsleahy} for more results along the same lines.

\smallskip
\newtheorem{siuthm}[qfcoh]{Theorem}
\begin{siuthm}
Let $X=(X, \ra)$ be a real-analytic manifold, and let $Y$ be a compact subset of $X.$ Then $\Gamma(Y, \ra)$ is Noetherian iff $Z \cap Y$ has only a finite number of topological components for every coherent real-analytic subvariety $Z$ of an open neighborhood $U$ of $Y$ in $X.$
\label{siutheorem}
\end{siuthm}

\smallskip As shown in \cite{lojasiewiczbook}, a compact semianalytic subset has only a finite number of topological components, so whenever $Y$ is a compact semianalytic subset of X, then $\Gamma(Y, \ra)$ is a Noetherian ring. We will now prove Theorem~\ref{finitetype} in separate propositions for $\tilde \II^q$ versus $\tilde \II^q_k$ for all $k \geq 1:$

\smallskip
\newtheorem{allcoh}[qfcoh]{Proposition}
\begin{allcoh}
Let $\Omega$ in $\C^n$ be a domain with real-analytic boundary $b \Omega.$ Let $\tilde U$ be any open subset of $b \Omega$ such that $\tilde U$ is contained in a compact semianalytic subset $Y$ of $b \Omega.$ If $\Omega$ is bounded, $b \Omega$ itself may be taken as $\tilde U.$ The ideal sheaf $\tilde \II^q$ of real-analytic subelliptic multipliers for the $\dbar$-Neumann problem on $(p,q)$ forms defined on $\tilde U$ is coherent. Viewed as a sheaf on all of $b \Omega,$ $\tilde \II^q$ is quasi-coherent.
\label{allcoherent}
\end{allcoh}

\smallskip\noindent {\bf Proof:} Define $\JJ^q= \tilde \II^q \otimes_{\ra} \E$ and seek to apply Tougeron's criterion from Proposition~\ref{qfcoherent}. To get quasi-flasqueness, use the argument in the proof of Proposition~\ref{bringhomeall}. $\JJ^q$ is of finite type by the observation following Theorem~\ref{siutheorem}. \qed

\smallskip
\newtheorem{stepscoh}[qfcoh]{Proposition}
\begin{stepscoh}
Let $\Omega$ in $\C^n$ be a pseudoconvex domain with real-analytic boundary $b \Omega.$ Let $\tilde U$ be any open subset of $b \Omega$ such that $\tilde U$ is contained in a compact semianalytic subset $Y$ of $b \Omega.$ If $\Omega$ is bounded, $b \Omega$ itself may be taken as $\tilde U.$ The multiplier ideal sheaf $\tilde \II^q_k$ given by the modified Kohn algorithm on $\tilde U$ at step $k$ for each $k \geq 1$ is coherent. Viewed as a sheaf on all of $b \Omega,$ $\tilde \II^q_k$ is quasi-coherent for each $k \geq 1.$
\label{stepscoherent}
\end{stepscoh}

\smallskip\noindent {\bf Proof:} We follow the same procedure here as in the proof of Proposition~\ref{allcoherent}. Define $\JJ^q_k= \tilde \II^q_k \otimes_{\ra} \E$ and apply Tougeron's criterion from Proposition~\ref{qfcoherent}. Quasi-flasqueness comes from the argument in the proof of Proposition~\ref{bringhomestepk}. $\JJ^q_k$ is of finite type by the observation following Theorem~\ref{siutheorem}. \qed

\bigskip\noindent {\bf Proof of Theorem~\ref{finitetype}} Combine the results in Propositions~\ref{allcoherent} and \ref{stepscoherent}. \qed

\bigskip The sheaf at the first step of the Kohn algorithm is the sheaf of all germs vanishing on the subvariety determined by all the Levi minors $coeff\{\partial r \wedge \dbar r \wedge (\partial \dbar r)^{n-q}\}$ inside the real-analytic manifold that is $b \Omega,$ so its support is simple to visualize. At subsequent steps, it is harder to visualize the support of $\tilde \II^q_k$ or for that matter, it is not at all clear how the support of $\tilde \II^q$ looks, especially if there are $q$-dimensional complex varieties in $b \Omega,$ so subellipticity fails. What is clear, however, is that these sheaves have better than expected properties because they are essentially coming from an estimate, namely estimate \eqref{subellest}. In Section 6 of \cite{kohnacta}, Kohn proved the following result, Proposition 6.5 on p.111, via standard techniques such as density of smooth points in a real-analytic variety, complexification, and Oka's coherence theorem:

\smallskip
\newtheorem{kohnseq}[qfcoh]{Proposition}
\begin{kohnseq}
If $I$ is an ideal of germs of real-analytic functions at $0 \in \R^p$ and if $I=\sqrt[\R]{I},$ then there exists a sequence of points $x^{(\nu)} \in \V(I)$ such that $x^{(\nu)}$ converges to $0$ and such that each $x^{(\nu)}$ has a neighborhood $U_\nu$ with the property that if $y \in U_\nu \cap \V(I),$ then $\I_y \V(I)$ is generated by the elements of $I.$
\label{kohnsequence}
\end{kohnseq}

\smallskip Kohn used this result in order to show the termination of his algorithm in the real-analytic case for a domain of finite D'Angelo type. Obviously, Theorem~\ref{finitetype} implies a much stronger result for $\II^q$ and $\tilde \II^q_k$ for each $k \geq 1$ as follows: 

\smallskip
\newtheorem{corft}[qfcoh]{Corollary}
\begin{corft}
Let $\Omega$ in $\C^n$ be a pseudoconvex domain with real-analytic boundary $b \Omega.$ For every $x_0 \in b \Omega,$ there exists a neighborhood $U_{x_0}$ around $x_0$ such that the sections $s_j$ that generate the stalk $\tilde \II^q(x_0)$ not only generate the stalk $\tilde \II^q(x)$ but also $\I \V(\tilde \II^q(x))$ for every $x \in U_{x_0}.$ The same is true for $\tilde \II^q_k$ for every $k \geq 1.$
\label{finitetypecorollary}
\end{corft}

\smallskip\noindent {\bf Proof:} By Theorem~\ref{finitetype}, the sheaves $\tilde \II^q$ and $\tilde \II^q_k$ for every $k \geq 1$ are quasi-coherent on $b \Omega,$ so each $x_0$ has a neighborhood $U_{x_0}$ where the generators of the stalk at $x_0$ generate the stalk at $x$ for every $x \in U_{x_0}.$ The assertion that those generators also generate $\I \V(\tilde \II^q(x))$ and $\I \V(\tilde \II^q_k(x))$ respectively is a consequence of the closure of $\tilde \II^q$ and $\tilde \II^q_k$ under the real radical (Theorem~\ref{subellpropcor} (b)) and of the \L ojasiewicz Nullstellensatz, which says that an ideal of real-analytic germs $I$ satisfies the Nullstellensatz, $I = \I \V (I),$ if $I=\sqrt[\R]{I};$ see \cite{lojasiewiczbook}. \qed

\bibliographystyle{plain}
\bibliography{SmoothTypeEquiv}

\end{document}